\documentclass[a4paper,11pt]{amsart}
\usepackage{amscd,amssymb, longtable}
\usepackage{amsthm,amsmath,amssymb}
\usepackage{supertabular}
\usepackage[matrix,arrow]{xy}
\usepackage{calc}
\usepackage{hyperref}

\theoremstyle{plain}
\newtheorem{theorem}[equation]{Theorem}

\newtheorem{proposition}[equation]{Proposition}
\newtheorem{lemma}[equation]{Lemma}

\theoremstyle{definition}

\newtheorem{definition}[equation]{Definition}
\newtheorem{problem}[equation]{Problem}

\theoremstyle{remark}

\renewcommand{\leq}{\leqslant}

\renewcommand{\geq}{\geqslant}

\theoremstyle{plain}
\newtheorem{innermytheorem}{Theorem}
\newenvironment{mytheorem}[1]
  {\renewcommand\theinnermytheorem{#1}\innermytheorem}
  {\endinnermytheorem}
\newtheorem{innermyproposition}{Proposition}
\newenvironment{myproposition}[1]
  {\renewcommand\theinnermyproposition{#1}\innermyproposition}
  {\endinnermyproposition}
\newtheorem{innermylemma}{Lemma}
\newenvironment{mylemma}[1]
  {\renewcommand\theinnermylemma{#1}\innermylemma}
  {\endinnermylemma}

\makeatletter
  \@addtoreset{equation}{section}
  \@addtoreset{section}{part}
\makeatother

\usepackage{listings}
\usepackage{color}

\definecolor{dkgreen}{rgb}{0,0.6,0}
\definecolor{gray}{rgb}{0.5,0.5,0.5}
\definecolor{mauve}{rgb}{0.58,0,0.82}

\DeclareRobustCommand{\lcm}{\operatorname{lcm}}
\DeclareRobustCommand{\Proj}{\operatorname{Proj}}
\DeclareRobustCommand{\Spec}{\operatorname{Spec}}
\DeclareRobustCommand{\wt}{\operatorname{wt}}
\DeclareRobustCommand{\GL}{\operatorname{GL}}

\newcommand{\condsai}{{\makebox[13.5pt]{$a_i$\hfill}}}
\newcommand{\condsaj}{{\makebox[13.5pt]{$a_j$\hfill}}}

\newcommand{\StateLemmaOne}{\textit{%
Let $(a_0, a_1, a_2, a_3, d)$ be an ordered tuple that satisfies conditions (i)--(iv) and $I = a_t + \frac{a_u}{2}$,
where one of $t$ and $u$ is equal to $2$ and the other less than $2$.
Then, either $I = a_i + a_j$ for some $i < j \leq 2$ or $I = a_0 + \frac{a_1}{2}$ or $I = \frac{a_0}{2} + a_1$.}
}

\newcommand{\StateLemmaTwo}{\textit{%
Let $(a_0, a_1, a_2, a_3, d)$ be an ordered tuple that satisfies conditions (i)--(iv) and $I = a_t + a_3$ or
$I = \frac{a_t}{2} + a_3$ or $I = a_t + \frac{a_3}{2}$, where $t \in \{0,1,2\}$.
Then, either $I = a_i + a_j$ for some $i < j \leq 2$ or $I = a_0 + \frac{a_1}{2}$ or $I = \frac{a_0}{2} + a_1$.}
}

\newcommand{\StateLemmaThree}{\textit{%
Let $(a_0, a_1, a_2, a_3, d)$ be an ordered tuple that satisfies conditions (i)--(iv) and $I = a_t + \frac{a_u}{2}$,
where $t,u \in \{0,1\}$ and $t \neq u$.
Then $a_3 = \frac{a_u}{2} + a_2$.}
}

\newcommand{\StatePropositionFourCondsSufficient}{\textit{%
Every tuple $(a_0, a_1, a_2, a_3, d)$ of type I, II or~III satisfying conditions (i)--(iv) also satisfies conditions (v) and~(vi),
that is, there exists a quasi-smooth well-formed hypersurface in $\mathbb{P}(a_0, a_1, a_2, a_3)$ of degree $d < a_0 + a_1 + a_2 + a_3$
and index $I = a_0 + a_1 + a_2 + a_3 - d$.}
}

\newcommand{\StateTheoremSeriesCheckOne}{\textit{%
If an ordered tuple $(b_0, b_1, b_2, b_3, e)$ is of type I, II or~III and satisfies conditions (i)--(iv),
then for every ordered tuple $(a_0, a_1, a_2, a_3, d)$ in its corresponding infinite series given in
Definition~\ref{definition:corresponding-series}, there exists a quasi-smooth well-formed hypersurface in
$\mathbb{P}(a_0, a_1, a_2, a_3)$ of degree $d < a_0 + a_1 + a_2 + a_3$ and index $I = a_0 + a_1 + a_2 + a_3 - d$}.
}

\author{Erik Paemurru}

\title{Del Pezzo Surfaces in Weighted Projective Spaces}

\thanks{The author was supported by a grant from the William Manson
bequest of the University of Edinburgh.}

\address{School of Mathematics, University of Edinburgh, Edinburgh, EH9 3JZ, UK}
\curraddr{Department of Mathematical Sciences, Loughborough University, Loughborough, LE11 3TU}
\email{E.Paemurru@lboro.ac.uk}

\begin{document}

\begin{abstract}
We study singular del Pezzo surfaces that are quasi-smooth and well-formed weighted hypersurfaces.
We give an algorithm how to classify all of them.
\end{abstract}

\maketitle

\section{Introduction} \label{section:intro}

Fano varieties are the building blocks of rationally connected varieties.
They have been studied for a long time and have often been used to produce counterexamples to long-standing conjectures.
Classically, Fano varieties were assumed to be smooth.
However, during the last decades the progress in the area gave tools and posed problems dealing with mildly singular Fano varieties.
The classification of singular Fano varieties seems to be absolutely hopeless in higher-dimensions, without bounding the singularities.
Nevertheless, we know many partial classification-type results about two-dimensional Fano varieties, also known as del Pezzo surfaces,
thanks to the combined effort of many people (see \cite{HidakaWatanabe}, \cite{MiyanishiZhang},
\cite{Zhang}, \cite{Kojima}, \cite{KeelMcKernan}, \cite{Shokurov}, \cite{AlNik}, \cite{Belousov}).

The classification problem for Fano manifolds is closely related
to the problem of the existence of K\"ahler--Einstein metrics on them (see \cite{Spotti}).
It has been conjectured by Yau, Tian and Donaldson that a Fano manifold admits a K\"ahler--Einstein metric if and only if it is $K$-polystable.
One direction of this conjecture, the $K$-polystability of K\"ahler--Einstein Fano manifolds,
follows from the works of Tian, Donaldson, Stoppa and Berman (see \cite{Ti97}, \cite{SD}, \cite{Stoppa}, \cite{Berman}).
The other direction has been recently proved by Chen, Donaldson and Sun in \cite{CDS1,CDS2,CDS3} and independently by Tian in \cite{Tian2015}.
Unfortunately, this result is not easy to apply since $K$-polystability is usually very hard to check.

The problem of the existence of a K\"ahler--Einstein metric on smooth del Pezzo surfaces has been explicitly solved by Tian in \cite{Ti90}.
For del Pezzo surfaces with quotient singularities, we do not have such an explicit solution (for orbifold metrics),
since del Pezzo surfaces with quotient singularities have not been classified.
So, it seems natural to consider this problem imposing some additional restrictions on the class of singular del Pezzo surfaces.
In this paper, we will consider singular del Pezzo surfaces that are quasi-smooth and well-formed
hypersurfaces in weighted projective spaces (see \cite[Definition~6.9]{IF00}).

Let $S_d$ be a hypersurface in $\mathbb{P}(a_0, a_1, a_2, a_3)$ of degree $d$,
where $1 \leq a_0 \leq a_1 \leq a_2 \leq a_3$ are some natural numbers.
Then $S_d$ is given by
\[
\phi(x, y, z, t) = 0 \subset \mathbb{P}(a_0, a_1, a_2, a_3) \cong \Proj(\mathbb{C}[x, y, z, t])
\]
where $\wt(x) = a_0$, $\wt(y) = a_1$, $\wt(z) = a_2$, $\wt(t) = a_3$,
and $\phi$ is a quasi-homogeneous polynomial of degree $d$ with respect to these weights.
The equation
\[
\phi(x, y, z, t) = 0 \subset\mathbb{C}^4 \cong \Spec(\mathbb{C}[x, y, z, t])
\]
defines a three-dimensional hypersurface quasi-homogeneous singularity $(V, O)$, where $O = (0, 0, 0, 0)$.
Recall that $S_d$ is called \emph{quasi-smooth} if the singularity $(V, O)$ is isolated.
Recall that $S_d$ is called \emph{well-formed} if
\[
\gcd(a_1, a_2, a_3) = \gcd(a_0, a_2, a_3) = \gcd(a_0, a_1, a_3) = \gcd(a_0, a_1, a_2) = 1
\]
and each positive integer $\gcd(a_0, a_1)$, $\gcd(a_0, a_2)$, $\gcd(a_0, a_3)$,
$\gcd(a_1, a_2)$, $\gcd(a_1, a_3)$, $\gcd(a_2, a_3)$ divides $d$.
If the hypersurface $S_d$ is quasi-smooth and well-formed, then it follows from \cite[Theorem~7.9]{Ko97}, \cite[Proposition~8.13]{Ko97},
\cite[Remark~8.14.1]{Ko97}, \cite[Theorem~11.1]{Ko97}, and the adjunction formula that the following conditions are equivalent
\begin{itemize}
\item the inequality $d < a_0 + a_1 + a_2 + a_3$ holds,
\item the singularity $(V, O)$ is a rational singularity,
\item the singularity $(V, O)$ is a Kawamata log terminal singularity,
\item the hypersurface $S_d$ is a del Pezzo surface with quotient singularities.
\end{itemize}

Starting from now, we set $d < \sum_{i=0}^n a_i$ and that the hypersurface $S_d$ is quasi-smooth and well-formed.
We define $I = a_0 + a_1 + a_2 + a_3 - d$.
This is usually called the index of the hypersurface $S_d \subset \mathbb{P}(a_0, a_1, a_2, a_3)$.
For every positive integer $I$ we have infinitely many possibilities for the sextuple $(a_0, a_1, a_2, a_3, d, I)$
such that there exists a quasi-smooth well-formed hypersurface in $\mathbb{P}(a_0, a_1, a_2, a_3)$ of degree $d < a_0 + a_1 + a_2 + a_3$
and index $I = a_0 + a_1 + a_2 + a_3 - d$.
This is not surprising, since we know there are infinitely many families of del Pezzo surfaces with quotient singularities.

\begin{problem}
\label{problem:classification} \emph{Describe} all sextuples $(a_0, a_1, a_2, a_3, d, I)$ such that there exists a quasi-smooth
well-formed hypersurface in $\mathbb{P}(a_0, a_1, a_2, a_3)$ of degree $d < a_0 + a_1 + a_2 + a_3$ and index $I = a_0 + a_1 + a_2 + a_3 - d$.
\end{problem}

This problem was posed by Orlov a long time ago in order to test his conjecture about
the existence of a full exceptional collection on del Pezzo surfaces with quotient singularities.
Later his conjecture was proved by Kawamata, Elagin, Ishii, and Ueda (see \cite{Kawamata}, \cite{Elagin}, \cite{IshiiUeda}).

The first step in solving Problem~\ref{problem:classification} was done by Johnson and Koll\'ar who proved the following:

\begin{theorem}[{\cite[Theorem~8]{JoKo01b}}] \label{theorem:Kollar-Johnson}
Suppose that $I=1$. Then
\begin{itemize}
\item either $(a_0, a_1, a_2, a_3, d) = (2, 2m+1, 2m+1, 4m+1, 8m+4)$ for some $m\in\mathbb{N}$,
\item or the quintuple $(a_0, a_1, a_2, a_3, d)$ lies in the sporadic set
\[
\left\{\aligned
&(1,1,1,1,3), (1,1,1,2,4), (1,1,2,3,6), (1,2,3,5,10), (1,3,5,7,15), (1,3,5,8,16),\\
&(2,3,5,9,18), (3,3,5,5,15), (3,5,7,11,25), (3,5,7,14, 28), (3,5,11,18, 36),\\
&(5,14,17,21,56), (5,19,27,31,81), (5,19,27,50,100), (7,11,27,37,81),\\
&(7,11,27,44,88), (9,15,17,20,60), (9,15,23,23,69), (11,29,39,49,127),\\
&(11,49,69,128,256), (13,23,35,57,127), (13,35,81,128,256)\\
\endaligned\right\}.
\]
\end{itemize}
Moreover, for each listed quintuple $(a_0, a_1, a_2, a_3, d)$, there exists a quasi-smooth well-formed hypersurface in
$\mathbb{P}(a_0, a_1, a_2, a_3)$ of degree $d = a_0 + a_1 + a_2 + a_3 - 1$.
\end{theorem}

The second step in solving Problem~\ref{problem:classification} was done by Cheltsov and Shramov who solved
Problem~\ref{problem:classification} for $I=2$ (see \cite[Corollary~1.13]{ChSh-Zoo}).

For Cheltsov, Johnson, Koll\'ar, and Shramov, the main motivation to prove Theorem~\ref{theorem:Kollar-Johnson} was
the Calabi problem for del Pezzo surfaces with quotient singularities and, in particular, the Calabi problem for
quasi-smooth well-formed hypersurfaces in $\mathbb{P}(a_0, a_1, a_2, a_3)$ of degree $d < a_0 + a_1 + a_2 + a_3$.
Regarding the latter, Gauntlett, Martelli, Sparks, and Yau proved

\begin{theorem}[\cite{GaMaSpaYau06}]
\label{theorem:Bishop-Lichnerowicz} The surface $S_d$ does not admit an orbifold K\"ahler--Einstein metric if $I > 3 a_0$.
\end{theorem}

Thus, the Calabi problem for quasi-smooth well-formed hypersurfaces in $\mathbb{P}(a_0, a_1, a_2, a_3)$
of degree $d < a_0 + a_1 + a_2 + a_3$ has a negative solution if $I > 3 a_0$.
On the other hand, Araujo, Boyer, Demailly, Galicki, Johnson, Koll\'ar, and Nakamaye proved that the Calabi problem for quasi-smooth
well-formed hypersurfaces in $\mathbb{P}(a_0, a_1, a_2, a_3)$ of degree $d = a_0 + a_1 + a_2 + a_3 - 1$ almost always has a positive solution.

\begin{theorem}[{\cite{DeKo01}, \cite{JoKo01b}, \cite{Ara02}, \cite{BoGaNa03}}] \label{theorem:BGN-CPS}
Suppose that $I = 1$. Then $S_d$ admits an orbifold K\"ahler--Einstein metric except possibly the case when
$(a_0, a_1, a_2, a_3, d) = (1, 3, 5, 7, 15)$ and the polynomial $\phi(x_0, x_1, x_2, x_3)$ does not contain the monomial $x_1 x_2 x_3$.
\end{theorem}

The proof of Theorem~\ref{theorem:BGN-CPS} implicitly uses the $\alpha$-invariant
introduced by Tian for smooth Fano varieties in \cite{Tian}.
For $S_d$, its algebraic counterpart can be defined as
\[
\alpha(S_d) = \sup \left\{ \lambda \in \mathbb{Q}\ \left|\ %
\begin{aligned}
& \text{the log pair $(S_d, \lambda D)$ is log canonical}\\
& \text{for every effective $\mathbb{Q}$-divisor $D\equiv-K_{S_d}$}\\
\end{aligned}
\right.\right\},
\]
and one can easily extend this definition to any Fano variety with at most Kawamata log terminal singularities.
Tian, Demailly, and Koll\'ar showed that the $\alpha$-invariant plays an important role in
the existence of an orbifold K\"ahler--Einstein metric on Fano varieties with quotient singularities
(see \cite{Tian}, \cite{DeKo01}, \cite{Ch07b}, \cite[Theorem~A.3]{ChSh08umn}).
In particular, we have

\begin{theorem}[{\cite{Tian}, \cite{DeKo01}, \cite[Theorem~A.3]{ChSh08umn}}] \label{theorem:KE}
If $\alpha(S_d) > 2/3$, then $S_d$ admits an orbifold K\"ahler--Einstein metric.
\end{theorem}

Araujo, Boyer, Demailly, Galicki, Johnson, Koll\'ar, and Nakamaye proved that $\alpha(S_d) > 2/3$ if $I = 1$ except exactly one case
when $(a_0, a_1, a_2, a_3, d) = (1, 3, 5, 7, 15)$ and the polynomial $\phi(x_0, x_1, x_2, x_3)$ does not contain
the monomial $x_1 x_2 x_3$ (in this case $\alpha(S_d) = 7/10 < 2/3$ by \cite{ChShPa08}).
A similar approach was used by Boyer, Cheltsov, Galicki, Nakamaye, Park, and Shramov for $I \geq 2$
(see \cite{BoGaNa03}, \cite{ChShPa08}, \cite{ChSh-Zoo}).

It is seems unlikely that Problem~\ref{problem:classification} has a \emph{nice} solution for all $I$ at once.
However, the results by Cheltsov, Johnson, Koll\'ar, and Shramov indicate it seems possible to solve
Problem~\ref{problem:classification} for any fixed $I$.
The main purpose of this paper is to prove this and to give an algorithm that solves
Problem~\ref{problem:classification} for any fixed $I$, that is,
finds the set of quintuples $(a_0, a_1, a_2, a_3, d)$ such that there exists
a quasi-smooth well-formed hypersurface in $\mathbb{P}(a_0, a_1, a_2, a_3)$ of degree $d < a_0 + a_1 + a_2 + a_3$
and index $I = a_0 + a_1 + a_2 + a_3 - d$.

We hope that our classification can be useful to produce a vast number of examples of non-K\"ahler--Einstein del Pezzo surfaces
with quotient singularities using different kind of existing obstructions.
For example, recently Spotti proved

\begin{theorem}[\cite{Spotti}]
\label{theorem:Spotti} Let $S$ be a del Pezzo surfaces with at most quotient singularities,
and let $N$ is the biggest natural number such that $S_d$ has a quotient singularity $\mathbb{C}^2 / G$ with $N = |G|$,
where $G$ is a finite subgroup in $\GL_2(\mathbb{C})$ that does not contain quasi-reflections.
Then $S$ does not admit an orbifold K\"ahler--Einstein metric if $K_S^2 N \geq 12$.
\end{theorem}

Using our classification, we immediately obtain a huge number of examples of quasi-smooth well-formed
hypersurfaces in $\mathbb{P}(a_0, a_1, a_2, a_3)$ of degree $d < a_0 + a_1 + a_2 + a_3$ that do not admit
an orbifold K\"ahler--Einstein metric by Theorem~\ref{theorem:Spotti} such that the obstruction found by
Gauntlett, Martelli, Sparks, and Yau, i.e.\ Theorem~\ref{theorem:Bishop-Lichnerowicz}, is not applicable.
The following tuples $(a_0, a_1, a_2, a_3, d)$ come from Table~\ref{table:ind-4}.
\begin{itemize}
\item The tuple $(1,3,4,8,12)$ from the series $(1, 3, 3a+1, 3b+2, 3a+3b+3)$ with $I = 4$ satisfies both
Theorem~\ref{theorem:Bishop-Lichnerowicz} and Spotti's inequality.
\item The tuple $(1,3,7,8,15)$ from the series $(1, 3, 3a+1, 3b+2, 3a+3b+3)$ with $I = 4$ satisfies
Theorem~\ref{theorem:Bishop-Lichnerowicz} but not Spotti's inequality.
\item The tuple $(2,2,3,7,10)$ from the series $(2, 2, 2a+1, 2b+1, 2a+2b+2)$ with $I = 4$ does not satisfy
Theorem~\ref{theorem:Bishop-Lichnerowicz} but satisfies Spotti's inequality.
\item The tuple $(2,2,3,3,6)$ from the series $(2, 2, 2a+1, 2b+1, 2a+2b+2)$ with $I = 4$ satisfies neither
Theorem~\ref{theorem:Bishop-Lichnerowicz} nor Spotti's inequality.
\end{itemize}
These examples show that the previous inequality by Gauntlett, Martelli, Sparks, and Yau (Theorem~\ref{theorem:Bishop-Lichnerowicz})
is independent of the new inequality discovered by Spotti (Theorem~\ref{theorem:Spotti}).
Thus, Spotti's inequality is a new and powerful obstruction to the existence of orbifold K\"ahler--Einstein metrics on
del Pezzo surfaces with quotient singularities.

Now, we state the main result of the paper.

\begin{theorem} \label{theorem:main-1}
The tuples $(a_0, a_1, a_2, a_3, d)$ such that there exists a quasi-smooth well-formed hypersurface in $\mathbb{P}(a_0, a_1, a_2, a_3)$
of degree $d < a_0 + a_1 + a_2 + a_3$ and index $I = a_0 + a_1 + a_2 + a_3 - d$ are precisely those that either
belong to Tables \ref{table:inf} and~\ref{table:spor} of the appendix, or are ordered,
satisfy conditions (i)--(iv) from Theorem~\ref{theorem:conds} and have one of the following forms where $x,y \in \mathbb{Z}$:
\begin{itemize}

\item class~1 tuples $(a_0, a_1, b_2 + xm, b_3 + ym, b_2 + b_3 + (x+y)m)$ where
\begin{itemize}
\item $a_0$, $a_1$ are positive integers such that $I = a_0 + a_1$,
\item $m = \lcm (a_0, a_1)$,
\item $b_2$, $b_3$ are positive integers,
\end{itemize}

\item class~2 tuples $(a_0, a_1, a_2, b_3 + xm, a_1 + b_3 + xm)$ where
\begin{itemize}
\item $a_0$, $a_1$, $a_2$ are positive integers such that $I = a_0 + a_2$ and $I > a_0 + a_1$,
\item $m = \lcm (a_0, a_1, a_2)$,
\item $b_3$ is a positive integer,
\end{itemize}

\item class~3 tuples $(a_0, a_1, a_2, b_3 + xm, a_0 + b_3 + xm)$ where
\begin{itemize}
\item $a_0$, $a_1$, $a_2$ are positive integers such that $I = a_1 + a_2$ and $I > a_0 + a_2$,
\item $m = \lcm (a_0, a_1, a_2)$,
\item $b_3$ is a positive integer,
\end{itemize}

\item class~4 tuples
$(a_0, a_1, b_2 + xm, \frac{a_1}{2} + b_2 + xm, a_1 + 2 b_2 + 2xm)$ where
\begin{itemize}
\item $a_0$, $a_1$ are positive integers such that $I = a_0 + \frac{a_1}2$,
\item $m = \lcm \left( a_0, \frac{a_1}2 \right)$,
\item $b_2$ is a positive integer,
\end{itemize}

\item class~5 tuples $(a_0, a_1, b_2 + xm, \frac{a_0}{2} + b_2 + xm, a_0 + 2 b_2 + 2xm)$ where
\begin{itemize}
\item $a_0$, $a_1$ are positive integers such that $I = \frac{a_0}2 + a_1$ and $I > a_0 + \frac{a_1}2$,
\item $m = \lcm \left( \frac{a_0}2, a_1 \right)$,
\item $b_2$ is a positive integer,
\end{itemize}

\item class~6 tuples $(a_0, a_1, b_2+xm, b_3+xm, d+2xm)$\\
$= (I-k, I+k, a+xm, a+k+xm, 2a+I+k+2xm)$, where
\begin{itemize}
\item $I$ is the index,
\item $k$ is a positive integer such that $I-k$ is positive,
\item $m = \lcm (a_0,a_1,k) = \lcm (I-k,I+k,k)$,
\item $a$ is a positive integer.
\end{itemize}

\end{itemize}
\end{theorem}

Quasi-smoothness and well-formedness conditions are given in Theorems \ref{theorem:basic-conds} and~\ref{theorem:conds}.
A tuple $(a_0, a_1, a_2, a_3, d)$ is ordered if $a_0 \leq a_1 \leq a_2 \leq a_3 \leq d$.
The fact that we only need to consider conditions (i) to (iv) is proved in Proposition~\ref{proposition:four-conds-sufficient}.

The second main result is the following:

\begin{theorem} \label{theorem:main-2}
If there exists a quasi-smooth well-formed hypersurface for a tuple from a class~$n$ series given in Theorem~\ref{theorem:main-1},
then there also exists a quasi-smooth well-formed hypersurface for all the other ordered tuples in the series, that is,
all the other ordered tuples with the same $a_i$ and $b_i$ but different $x$ or~$y$.
\end{theorem}

This result is proved in Theorem~\ref{theorem:series-check-one}.

For any fixed $I$, there are only finitely many series in Theorem~\ref{theorem:main-1}.
Let us check this for class~1.
A series in class~1 is given by ordered tuples $(a_0, a_1, b_2 + xm, b_3 + ym, b_2 + b_3 + (x+y)m)$
where $x,y \in \mathbb{Z}$ and $a_0, a_1, b_2, b_3, m$ are fixed positive integers.
Since $a_0 + a_1 = I$, there are only finitely many choices for $a_0$ and $a_1$.
The number $m = \lcm(a_0,a_1)$ is uniquely determined by $a_0$ and $a_1$.
Since $x,y \in \mathbb{Z}$, we only care for $b_2$ and $b_3$ modulo $m$.
A series is uniquely determined by $a_0$, $a_1$, $b_2$, $b_3$ and $m$, so there are only finitely many series in class~1.
Similarly for classes $2$ to~$6$.

According to Theorem~\ref{theorem:main-2}, checking one tuple from every such series determines
whether every tuple in the series is such that there exists a quasi-smooth well-formed hypersurface.
So, the algorithm to classify the hypersurfaces for an index $I$ reduces to checking
conditions (i)--(iv) of Theorem~\ref{theorem:conds} for finitely many tuples.

In Section~\ref{section:main}, we prove Theorems \ref{theorem:main-1} and~\ref{theorem:main-2},
leaving computations to Section~\ref{section:proofs}.
The main tool we use is Theorem~\ref{theorem:conds} which lists the quasi-smoothness and well-formedness conditions.
The complete lists of the classified hypersurfaces for indices $I$ from $1$ to~$6$ are given in Section~\ref{section:small-index-cases}.
Appendix~\ref{section:source-code} contains the computer code of the algorithm.

Theorems \ref{theorem:main-1} and~\ref{theorem:main-2} provide the algorithm to solve Orlov's problem,
Problem~\ref{problem:classification}, for any fixed $I$, as well as giving the general form of the answer for any $I$.
The surprisingly rigid form describes the hypersurfaces for all indices at once, without needing to calculate them explicitly.

\bigskip
\textbf{Acknowledgements.}
I would like to thank the William Manson bequest which supported the project financially.

\bigskip
\section{Technical result} \label{section:main}

The following theorem describes the tuples $(a_0, a_1, a_2, a_3, d)$ such that there exists a quasi-smooth well-formed hypersurface.
Conditions (iv), (v) and~(vi) are taken from \cite[Conditions 2.1, 2.2 and 2.3]{JoKo01b}.

\begin{theorem} \label{theorem:basic-conds}
There exists a quasi-smooth well-formed hypersurface in $\mathbb{P}(a_0, a_1, a_2, a_3)$ of degree
$d < a_0 + a_1 + a_2 + a_3$ and index $I = a_0 + a_1 + a_2 + a_3 - d$ iff all of the following conditions hold:
\begin{itemize}
\item[(i--ii)] the hypersurface is well-formed,
\item[(iii)] the hypersurface is not degenerate, that is, none of the weights $a_i$ equal the degree of the hypersurface,
\item[(iv)] for every $i$, there exists $j$ ($j$ may equal $i$) such that there exists a monomial $x_i^{m_i} x_j$ of degree $d$
where $m_i \geq 1$,
\item[(v)] for every $i < j$ such that $\gcd (a_i, a_j) > 1$, there is a monomial $x_i^{b_i}x_j^{b_j}$ of degree $d$
where $b_i,b_j \geq 0$ and $b_i + b_j \geq 2$,
\item[(vi)] for every $i < j$, either there is a monomial $x_i^{b_i} x_j^{b_j}$ of degree $d$ or there exist $k < l$
such that the indices $i,j,k,l$ are pairwise different and there are monomials $x_i^{c_i}x_j^{c_j}x_k$ and $x_i^{d_i}x_j^{d_j}x_l$
of degree $d$ where $b_i, b_j, c_i, c_j, d_i, d_j \geq 0$ and $b_i + b_j \geq 2$ and $c_i + c_j \geq 1$ and $d_i + d_j \geq 1$.
\end{itemize}
\end{theorem}

We will use the following form of the theorem, adding results from \cite[Theorem 2.3]{ChSh-Zoo} and \cite[Definitions 1.10 and 2.2]{ChSh-Zoo}.

\begin{theorem} \label{theorem:conds}
There exists a quasi-smooth well-formed hypersurface in $\mathbb{P}(a_0, a_1, a_2, a_3)$ of degree $d < a_0 + a_1 + a_2 + a_3$
and index $I = a_0 + a_1 + a_2 + a_3 - d$ iff all of the following conditions hold:
\begin{itemize}
\item[(i)] for every $i < j$, $\gcd(a_i, a_j) \mid d$,
\item[(ii)] for every $i < j < k$, $\gcd (a_i, a_j, a_k) = 1$,
\item[(iii)] $d > a_3$,
\item[(iv)] for every $i$, there exists $j$ ($j$ may equal $i$) such that $a_i \mid d - a_j$,

\item[(v)] for every $i < j$ such that $\gcd (a_i,a_j) > 1$, one of the following holds:
\begin{itemize}
\item \condsai $\mid d$,
\item \condsaj $\mid d$,
\item \condsai $\mid d - a_j$,
\item \condsaj $\mid d - a_i$,
\item there exists $b_j \geq 2$ such that $a_i \mid d - a_j b_j$ and $d - a_j b_j \geq 0$,
\end{itemize}

\item[(vi)] for every $i < j$, (at least) one of the following holds:
\begin{itemize}
\item one of the following holds:
\begin{itemize}
\item \condsai $\mid d$,
\item \condsaj $\mid d$,
\item \condsai $\mid d - a_j$,
\item \condsaj $\mid d - a_i$,
\item there exists $b_j \geq 2$ such that $a_i \mid d - a_j b_j$ and $d - a_j b_j \geq 0$,
\end{itemize}
\item for pairwise different indices $i,j,k,l$ satisfying $k<l$, both of the following hold:
\begin{itemize}

\item one of the following holds:
\begin{itemize}
\item \condsai $\mid d - a_k$,
\item \condsaj $\mid d - a_k$,
\item there exists $c_j \geq 1$ such that $a_i \mid d - a_k - a_j c_j$ and $d - a_k - a_j c_j \geq 0$,
\end{itemize}

\item one of the following holds:
\begin{itemize}
\item \condsai $\mid d - a_l$,
\item \condsaj $\mid d - a_l$,
\item there exists $d_j \geq 1$ such that $a_i \mid d - a_l - a_j d_j$ and $d - a_l - a_j d_j \geq 0$,
\end{itemize}

\end{itemize}
\end{itemize}
\end{itemize}

and one of the following conditions holds:
\begin{itemize}
\item[(type~I)] $I = a_i + a_j$ for some $i \neq j$,
\item[(type~II)] $I = a_i + \frac{a_j}{2}$ for some $i \neq j$,
\item[(type~III)] $(a_0, a_1, a_2, a_3, d) = (I-k, I+k, a, a+k, 2a+I+k)$, where $1 \leq k < I$ and $a \geq I + k$,
\item[(type~IV)] $(a_0, a_1, a_2, a_3, d)$ belongs to one of the infinite series listed in Table~\ref{table:inf}
or to the sporadic set given in Table~\ref{table:spor} of the appendix.
\end{itemize}
\end{theorem}

Throughout the paper, whenever refer to condition~(x), we refer to condition~(x) of Theorem~\ref{theorem:conds}.
Whenever we refer to type~X, we refer to type~X of Theorem~\ref{theorem:conds}.
Tuples of type~IV are simply read from the tables, so the rest of the paper concerns with types I--III.

Type~III of Theorem~\ref{theorem:conds} was changed from $0 \leq k < I$ in~\cite[Theorem 2.2]{ChSh-Zoo} to $1 \leq k < I$,
as the only tuple from the case $k = 0$ already exists in the tables for type IV.

The aim of the paper is to find all the tuples $(a_0, a_1, a_2, a_3, d)$
such that there exists a quasi-smooth well-formed hypersurface, given a fixed index $I$.
At first, it is easier to deal with a superset of such tuples, namely, when only conditions (i)--(iv) are satisfied.
Proposition~\ref{proposition:four-conds-sufficient} shows that for tuples of type I--III this weaker set of conditions suffices,
that is, for every such tuple there exists a quasi-smooth well-formed hypersurface.

Below, we consider a certain subset of the tuples of types I--III.
If a tuple belongs to this subset, we assign a unique \emph{class number} and an infinite series to it.

\begin{definition} \label{definition:classes}
Given an index $I$, an ordered tuple $(a_0, a_1, a_2, a_3, d)$ is assigned a unique \emph{class number} if it satisfies one of the following:

\begin{itemize}
\item class~1 --- $I = a_0 + a_1$,
\item class~2 --- $I = a_0 + a_2$ and $I > a_0 + a_1$,
\item class~3 --- $I = a_1 + a_2$ and $I > a_0 + a_2$,
\item class~4 --- $I = a_0 + \frac{a_1}{2}$ and $a_3 = \frac{a_1}2 + a_2$,
\item class~5 --- $I = \frac{a_0}{2} + a_1$ and $I > a_0 + \frac{a_1}{2}$ and $a_3 = \frac{a_0}2 + a_2$,
\item class~6 --- The tuple is of type~III, that is, it satisfies\\
$(a_0, a_1, a_2, a_3, d) = (I-k, I+k, a, a+k, 2a+I+k)$,\\
where $1 \leq k < I$ and $a \geq I + k$.
\end{itemize}
\end{definition}

\begin{definition} \label{definition:corresponding-series}
Infinite series for tuples of class 1--6.
\begin{itemize}
\item For class~1 tuples, the corresponding series is $(a_0, a_1, a_2 + xm, a_3 + ym, a_2 + a_3 + (x+y)m)$,
where $m = \lcm (a_0, a_1)$ and $x,y \in \mathbb{Z}$.
\item For class~2 tuples, the corresponding series is $(a_0, a_1, a_2, a_3 + xm, a_1 + a_3 + xm)$,
where $m = \lcm (a_0, a_1, a_2)$ and $x \in \mathbb{Z}$.
\item For class~3 tuples, the corresponding series is $(a_0, a_1, a_2, a_3 + xm, a_0 + a_3 + xm)$,
where $m = \lcm (a_0, a_1, a_2)$ and $x \in \mathbb{Z}$.
\item For class~4 tuples, the corresponding series is $(a_0, a_1, a_2 + xm, \frac{a_1}{2} + a_2 + xm, a_1 + 2 a_2 + 2xm)$,
where $m = \lcm \left( a_0, \frac{a_1}2 \right)$ and $x \in \mathbb{Z}$.
\item For class~5 tuples, the corresponding series is $(a_0, a_1, a_2 + xm, \frac{a_0}{2} + a_2 + xm, a_0 + 2 a_2 + 2xm)$,
where $m = \lcm \left( \frac{a_0}2, a_1 \right)$ and $x \in \mathbb{Z}$.
\item For class~6 tuples, the corresponding series is $(a_0, a_1, a_2+xm, a_3+xm, d+2xm) = (I-k, I+k, a+xm, a+k+xm, 2a+I+k+2xm)$,
where $m = \lcm (a_0,a_1,k) = \lcm (I-k,I+k,k)$ and $x \in \mathbb{Z}$.
\end{itemize}
The parameters $x$ and $y$ are bounded below such that the tuple is ordered.
\end{definition}

By an infinite series $(a_0, a_1, a_2 + xm, a_3 + ym, a_2 + a_3 + (x+y)m)$ for a class~1 tuple,
we mean the set $\left\{ (a_0, a_1, a_2 + xm, a_3 + ym, a_2 + a_3 + (x+y)m) \mid x, y \in \mathbb{Z} \right\}$
where all the tuples $(b_0, b_1, b_2, b_3, e)$ in the set are ordered, that is, $b_0 \leq b_1 \leq b_2 \leq b_3 \leq e$.
Analogously for classes 2--6.

Below, we prove Proposition~\ref{proposition:conds-imply-class} which implies that every tuple of type I--III
for which there exists a quasi-smooth well-formed hypersurface is of some class 1--6.
In fact, we show more: every tuple of type I--III satisfying conditions (i)--(iv) is of some class 1--6.
First, we state three lemmas, the proofs of which are in Section~\ref{section:proofs}.

\begin{lemma} \label{lemma:lemma-1}
\StateLemmaOne
\end{lemma}

\begin{lemma} \label{lemma:lemma-2}
\StateLemmaTwo
\end{lemma}

\begin{lemma} \label{lemma:lemma-3}
\StateLemmaThree
\end{lemma}

\begin{proposition} \label{proposition:conds-imply-class}
An ordered tuple of type I, II or III satisfying conditions (i)--(iv) is of some class 1--6.

\begin{proof}
Follows directly from the three lemmas above.
\end{proof}
\end{proposition}

Using Proposition~\ref{proposition:conds-imply-class}, we can show the
following.

\begin{proposition} \label{proposition:four-conds-sufficient}
\StatePropositionFourCondsSufficient
\end{proposition}

Using the last two propositions, we can prove the main theorem.

\begin{theorem} \label{theorem:series-check-one}
\StateTheoremSeriesCheckOne
\end{theorem}

The proofs of Proposition~\ref{proposition:four-conds-sufficient}
and Theorem~\ref{theorem:series-check-one} are in Section~\ref{section:proofs}.

Proposition~\ref{proposition:conds-imply-class} implies that all the tuples for which there exists a quasi-smooth well-formed hypersurface
lie in the infinite series given in Definition~\ref{definition:corresponding-series} or belong to Tables \ref{table:inf} or~\ref{table:spor}.
Theorem~\ref{theorem:series-check-one} says that we only need to check one tuple for conditions (i)--(iv) to see
whether all the tuples in the series are such that there exists a quasi-smooth well-formed hypersurface.
As there are only finitely many infinite series for any index $I$ (by the argument given below Theorem~\ref{theorem:main-2}),
this provides an effective algorithm to classify the hypersurfaces for a fixed $I$.

\bigskip
\section{Proofs} \label{section:proofs}

To prove Lemmas \ref{lemma:lemma-1} and~\ref{lemma:lemma-2}, we use condition~(iv) of Theorem~\ref{theorem:conds}
choosing $i = 3$ to express the weight $a_3$ in terms of $a_0$, $a_1$ and~$a_2$.
Then, we try to express the index $I = a_0 + a_1 + a_2 + a_3 - d$ in terms of the weights as in the statement of the lemma.
If unsuccessful, we use condition~(iv) again, choosing $i = 2$ to express $a_2$ in terms of $a_0$ and~$a_1$.
The weights $a_0, a_1, a_2, a_3$ are ordered integers with $1 \leq a_0 \leq a_1 \leq a_2 \leq a_3$.

\begin{mylemma}{\ref*{lemma:lemma-1}} \label{proof:lemma-1}
\StateLemmaOne

\begin{proof}
Let us define the index $v$ with
\begin{align*}
\{t,u,v\} & = \{0,1,2\}\\
v & \in \{0,1\}.
\end{align*}
We find
\[
d = a_v + \frac{a_u}{2} + a_3.
\]
Using condition~(iv) with $i = 3$, we find $a_3 \mid d - a_j$, giving $a_3 \mid a_v + \frac{a_u}{2} - a_j$ where $j \in \{t,u,v,3\}$.
Noting that $-a_3 < a_v + \frac{a_u}{2} - a_j < 2 a_3$, we can define $x$ such that
\begin{align*}
x & \in \{0,a_3\}\\
x & = a_v + \frac{a_u}{2} - a_j.
\end{align*}
The cases $(j,x) \in \{(v,0), (v,a_3), (u,a_3), (3,a_3)\}$ easily give contradictions.
If $(j,x) = (u,0)$, then $\frac{a_u}{2} = a_v$ and $I = a_t + a_v$.
If $(j,x) = (t,0)$, then $a_t = a_v + \frac{a_u}{2}$ and $I = a_u + a_v$.
It is left to consider $(j,x) \in \{(3,0), (t,a_3)\}$.
For both of these, we can define $y$ such that
\begin{align*}
y & \in \{0, a_t\}\\
a_3 & = \frac{a_u}{2} + a_v - y\\
d & =  a_u + 2 a_v - y\\
I & = a_t + \frac{a_u}{2}.
\end{align*}
Using condition~(iv) with $i = 2$, we find $a_2 \mid a_u + 2 a_v - y - a_k$ where $k \in \{t,u,v,3\}$.
It is easy to see that $a_u + 2 a_v - y - a_k \in \{a_2, 2 a_2\}$.
We consider these cases separately.
\smallskip

First, we consider $a_1 = a_2$.
If $\{t, u\} = \{0, 2\}$, then $I = a_0 + \frac{a_1}{2}$ or $I = \frac{a_0}{2} + a_1$.
This leaves the case $\{t, u\} = \{1, 2\}$.
We find
\begin{align*}
y & \in \{0, a_1\}\\
a_3 & = \frac{a_1}{2} + a_0 - y\\
d & =  a_1 + 2 a_0 - y\\
I & = \frac{3}{2} a_1.
\end{align*}
Since $I = \frac{3}{2} a_1$, we find $2 \mid a_1$.
Using condition~(i) and $a_1 = a_2$, we find $a_1 \mid 2 a_0$.
Since $\frac{a_1}2 \mid a_0, a_1, a_2$, condition~(ii) implies $(a_0, a_1, a_2) = (1, 2, 2)$.
We find $a_3 = 2$, and condition~(ii) gives us a contradiction.
\smallskip

Next, we consider $a_u + 2 a_v - y - a_k = 2 a_2$ and $a_1 < a_2$.
By definition $v \in \{0,1\}$.
Since $v = 0$ gives a contradiction, we find $v = 1$.
Similarly, $k \geq 1$ gives a contradiction.
We find $(k,t,v,u) = (0,0,1,2)$, giving
\begin{align*}
y & \in \{0, a_0\}\\
a_2 & = 2 a_1 - a_0 - y\\
a_3 & = - \frac{a_0}{2} + 2 a_1 - \frac{3}{2}y\\
d & = - a_0 + 4 a_1 - 2 y\\
I & = \frac{a_0}{2} + a_1 - \frac{y}{2}.
\end{align*}
If $y = 0$, then $I = \frac{a_0}2 + a_1$.
This leaves the case $y = a_0$.
We find $a_3 = 2 a_1 - 2 a_0 = a_2$ and $d = 4 a_1 - 3 a_0$.
Condition~(i) implies $a_2 \mid a_0$, giving $a_1 = a_2$, a contradiction.
\smallskip

Finally, we consider $a_u + 2 a_v - y - a_k = a_2$ and $a_1 < a_2$.
We find $2 a_3 + y = a_2 + a_k$, giving $y = 0$ and $a_k = a_2 = a_3$.
Therefore
\begin{align*}
a_3 & = a_2 = \frac{a_u}2 + a_v\\
d & =  a_u + 2 a_v\\
I & = a_t + \frac{a_u}{2}.
\end{align*}
If $t = 2$, then $I = a_u + a_v$.
This leaves the case $u = 2$.
We find $a_2 = 2 a_v$.
Since $a_v \mid a_v, a_2, a_3$, condition~(ii) implies $(a_0, a_1, a_2, a_3) = (1,1,2,2)$.
We find $I = a_0 + a_1$.
\end{proof}
\end{mylemma}

\begin{mylemma}{\ref*{lemma:lemma-2}} \label{proof:lemma-2}
\StateLemmaTwo

\begin{proof}
We define $u,v$ with
\begin{align*}
\{t,u,v\} & = \{0,1,2\}\\
u & < v.
\end{align*}

First, let us consider $I = a_t + a_3$.
We find
\[
d = a_u + a_v.
\]
Condition~(iv) gives $a_3 \mid d - a_j$ where $j \in \{0,1,2,3\}$.
Condition~(iii) gives $d > a_3$, and we find
\[
d - a_j = a_3.
\]
This implies $I = a_t + a_u + a_v - a_j$.
If $j \in \{0,1,2\}$, we find $I = a_p + a_q$ where $p < q \leq 2$, as required.
This leaves the case $j = 3$, giving $a_u + a_v = 2 a_3$.
We find $a_u = a_v = a_3$ and condition~(ii) implies $(a_0, a_1, a_2, a_3) = (1, 1, 1, 1)$.
So, $I = a_0 + a_1$.
\smallskip

Next, let us consider $I = \frac{a_t}{2} + a_3$.
We find
\[
d = a_u + a_v + \frac{a_t}2.
\]
Condition~(iv) gives us $a_3 \mid a_u + a_v + \frac{a_t}2 - a_j$ where $j \in \{t,u,v,3\}$.
We find
\[
a_3 = a_u + a_v + \frac{a_t}2 - a_j.
\]
This implies $I = a_t + a_u + a_v - a_j$.
So, it suffices to consider $j = 3$, giving
\[
d = 2 a_3.
\]
Using condition~(iv), we find $a_2 \mid a_u + a_v + \frac{a_t}2 - a_k$ where $k \in \{t,u,v,3\}$, giving $a_2 = a_u + a_v + \frac{a_t}2 - a_k$.
Since $a_2 \leq a_3$, we find $a_2 = a_3$, giving $d = a_2 + a_3$.
From the definition of the index, we find $I = a_0 + a_1$.
\smallskip

Finally, we consider $I = a_t + \frac{a_3}2$.
We find
\[
d = a_u + a_v + \frac{a_3}2.
\]
Condition~(iv) gives $a_3 \mid a_u + a_v + \frac{a_3}2 - a_j$ where $j \in \{t,u,v,3\}$.
So, we can define $x$ such that
\begin{align*}
x & \in \left\{ \frac{a_3}2, \frac{3 a_3}2 \right\}\\
x & = a_u + a_v - a_j.
\end{align*}
If $(j,x) \in \left\{ \left(t,\frac{a_3}2\right), \left(u,\frac{a_3}2\right), \left(v,\frac{a_3}2\right) \right\}$,
then $I = a_p + a_q$ where $p < q \leq 2$, as required.
The cases $(j,x) \in \left\{ \left(u,\frac{3 a_3}2\right), \left(v,\frac{3 a_3}2\right),
\left(3,\frac{3 a_3}2\right)\right\}$ give contradictions.
This leaves the cases $(j,x) \in \left\{ \left(3,\frac{a_3}2\right), \left(t,\frac{3 a_3}2\right) \right\}$.
For both of these, we can define $y$ such that
\begin{align*}
y & \in \{0, a_t\}\\
a_3 & = \frac{2 a_u + 2 a_v - 2 y}{3}\\
d & = \frac{4 a_u + 4 a_v - y}{3}.
\end{align*}
Condition~(iv) implies $a_2 \mid \frac{4 a_u + 4 a_v - y - 3 a_k}3$ where $k \in \{0,1,2,3\}$.
We can define $z$ such that
\begin{align*}
z & \in \{a_2, 2 a_2\}\\
z & = \frac{4 a_u + 4 a_v - y - 3 a_k}{3}.
\end{align*}
We consider the values of $z$ separately.

\begin{itemize}
\item If $z = a_2$, then $d = 2 a_3 + y = a_k + a_2$. We find $y = 0$ and $a_2 = a_3$, giving $d = a_2 + a_3$.
From the definition of the index, we find $I = a_0 + a_1$.
\item If $z = 2 a_2$, then $a_2 \geq a_u, a_v, a_t$ implies $k = 0$ and $k \neq u$.
By definition $u < v$, so $(k,t,u,v) = (0,0,1,2)$.
We find
\begin{align*}
y & \in \{0, a_0\}\\
a_2 & = \frac{- 3 a_0 + 4 a_1 - y}2\\
a_3 & = - a_0 + 2 a_1 - y\\
d & = - 2 a_0 + 4 a_1 - y\\
I & = \frac{a_0}{2} + a_1 - \frac{y}2.
\end{align*}
If $y = 0$, then $I = \frac{a_0}2 + a_1$.
If $y = a_0$, we find
\begin{align*}
a_2 & = a_3 = 2 a_1 - 2 a_0\\
d & = 4 a_1 - 3 a_0.
\end{align*}
Condition~(i) implies $a_3 \mid a_0$ and condition~(ii) implies $a_3 = 1$, a contradiction, since $2 \mid a_3$.
\end{itemize}
\end{proof}
\end{mylemma}

\begin{mylemma}{\ref*{lemma:lemma-3}} \label{proof:lemma-3}
\StateLemmaThree

\begin{proof}
We have
\[
d = \frac{a_u}{2} + a_2 + a_3.
\]
Using condition~(iv), we find $a_3 \mid \frac{a_u}{2} + a_2 - a_j$ where $j \in \{0,1,2,3\}$, giving $\frac{a_u}{2} + a_2 - a_j \in \{0, a_3\}$.
If $\frac{a_u}{2} + a_2 - a_j = 0$, we find $a_3 = \frac{a_u}{2} + a_2$ as required.
This leaves the case $\frac{a_u}{2} + a_2 - a_j = a_3$, giving $u = 1$ and $j = 0$.
We have
\begin{align*}
a_3 & = -a_0 + \frac{a_1}{2} + a_2\\
d & = -a_0 + a_1 + 2 a_2.
\end{align*}
Using condition~(iv), we find $a_2 \mid -a_0 + a_1 - a_k$ where $k \in \{0,1,2,3\}$, giving $-a_0 + a_1 - a_k \in \{-a_2, 0\}$.
Since $a_3 \geq a_k$, the case $-a_0 + a_1 - a_k = - a_2$ gives a contradiction.
Therefore, we find $a_1 = a_0 + a_k$, giving
\begin{align*}
a_1 & = 2 a_0\\
a_3 & = a_2\\
d & = a_0 + 2 a_2.
\end{align*}
Condition~(i) implies $a_2 \mid a_0$.
Condition~(ii) implies $a_2 = 1$, a contradiction, since $2 \mid a_1$.
\end{proof}
\end{mylemma}

\begin{myproposition}{\ref*{proposition:four-conds-sufficient}} \label{proof:four-conds-sufficient}
\StatePropositionFourCondsSufficient

\begin{proof}
By Proposition~\ref{proposition:conds-imply-class}, it suffices to consider tuples of class 1--6.
\smallskip

\emph{Classes 1--3}: $I = a_t + a_u$, where $t < u \leq 2$.
We can define $v$ such that
\begin{align*}
\{t,u,v\} & = \{0,1,2\}\\
d & = a_v + a_3.
\end{align*}

If $(i,j) = (v,3)$, then conditions (v) and~(vi) are satisfied, since $a_3 \mid d - a_v$.

If $(i,j) = (t,u)$, then condition~(iv) implies $a_u \mid d - a_k$ where $k \in \{t,u,v,3\}$.
If $k \in \{t,u\}$, conditions (v) and~(vi) hold.
If $k \in \{v,3\}$, then $a_u \mid a_v$ or $a_u \mid a_3$.
Condition~(i) implies $a_u \mid a_v, a_3$ and (ii) implies $a_u = 1$.
Therefore $a_u \mid d$, so (v) and~(vi) hold.

If $(i,j) \neq (t,u), (v,3)$, then either $\gcd(a_i, a_j) \mid a_v$ or $\gcd(a_i, a_j) \mid a_3$.
Condition~(i) implies $\gcd(a_i, a_j) \mid a_v, a_3$, and (ii) implies $\gcd(a_i, a_j) = 1$, so (v) is satisfied.

It is left to consider condition~(vi) for pairs $(i,j) \neq (t,u), (v,3)$.
We note the order of $i,j$ is not important in~(vi).
Similarly, the order of $k,l$ is not important in the second part of~(vi).
So, it suffices to consider
\begin{align*}
i & \in \{t,u\}\\
j & \in \{v,3\}.
\end{align*}
Using condition~(iv), we find
\[
a_i \mid d - a_k
\]
where $k \in \{0,1,2,3\}$.
If $k \in \{i,j\}$, then (vi) is satisfied.
If $k \in \{v,3\}$, then $a_i \mid a_v$ or $a_i \mid a_3$, and as before we find $a_i \mid d$, so (vi) is satisfied.
This leaves the case
\begin{align*}
k & \in \{t,u\}\\
k & \neq i.
\end{align*}
We define $l$ with
\begin{align*}
l & \in \{v,3\}\\
l & \neq j.
\end{align*}
The indices $i,j,k,l$ are pairwise different and $a_i \mid d - a_k$ and $a_j \mid d - a_l$.
So, by the second part of~(vi), condition~(vi) holds.
\smallskip

\emph{Classes 4--5}: $I = a_t + \frac{a_u}{2}$, where $t,u \in \{0,1\}$ with $t \neq u$.
By Lemma~\ref{lemma:lemma-3}, we have
\begin{align*}
a_3 & = \frac{a_u}{2} + a_2\\
d & = 2 a_3 = a_u + 2 a_2.
\end{align*}

If $j = 3$, then $a_j \mid d$, so conditions (v) and~(vi) are satisfied.

If $(i,j) = (u,2)$, then $a_2 \mid d - a_u$, so (v) and~(vi) hold.

If $(i,j) = (t,2)$, then (iv) implies $a_t \mid d - a_p$ where $p \in \{t,u,2,3\}$.
If $p \in \{t,2,3\}$, then (v) and~(vi) are satisfied.
If $p = u$, then $a_t \mid 2 a_2$.
Using conditions (i) and~(ii), we find $a_t \in \{1,2\}$, giving $a_t \mid d$.
So, (v) and~(vi) are satisfied.

If $(i,j) = (t,u)$, ignoring the order of $i$ and $j$, then (iv) implies $a_u \mid d - a_q$ where $q \in \{t,u,2,3\}$.
If $q \in \{t,u,3\}$, then (v) and~(vi) are satisfied.
If $q = 2$, then $a_u \mid a_2$, giving $a_u \mid d$, so (v) and~(vi) are satisfied.
\smallskip

\emph{Class 6}: $(a_0, a_1, a_2, a_3, d) = (I-k, I+k, a, a+k, 2a+I+k)$ where $1 \leq k < I$ and $a \geq I + k$.
We have
\[
d = a_0 + 2 a_3 = a_1 + 2 a_2.
\]

If $(i,j) \in \{(0,3), (1,2)\}$, then (v) and~(vi) satisfied.

If $(i,j) = (2,3)$, we find $a_3 \mid d - a_0$ and $a_2 \mid d - a_1$.
By denoting $(i,j,k,l) = (2,3,0,1)$ in the second part of~(vi), condition~(vi) holds.
For condition~(v), we note that (i) implies $\gcd(a_2,a_3) \mid a_0, a_1$, therefore $\gcd(a_2,a_3) = 1$ and (v) is satisfied.

If $(i,j) \in \{(0,1), (0,2), (1,3)\}$, then (iv) implies
\[
a_0 \mid d - a_p
\]
where $p \in \{0,1,2,3\}$. We show that for any $p$, we have either $a_0 \mid d$ or $a_0 \mid d - a_2$.
\begin{itemize}
\item If $p \in \{0,3\}$, then $a_0 \mid d$.
\item If $p = 1$, then (i) and~(ii) imply $a_0 \in \{1,2\}$, giving $a_0 \mid d$.
\item If $p = 2$, then $a_0 \mid d - a_2$.
\end{itemize}
Similarly, we can show that either $a_1 \mid d$ or $a_1 \mid d - a_3$.
Therefore, conditions (v) and~(vi) are satisfied for $(i,j) \in \{(0,2), (1,3)\}$.
Also, conditions (v) and~(vi) are satisfied for $(i,j) = (0,1)$ if either $a_0 \mid d$ or $a_1 \mid d$.

It is left to consider $(i,j) = (0,1)$ with
\begin{align*}
a_0 & \mid d - a_2\\
a_1 & \mid d - a_3.
\end{align*}
Choosing $(i,j,k,l) = (0,1,2,3)$, we see by the second part of condition~(vi) that (vi) is satisfied.
Next, we check~(v). From above, we see that
\begin{align*}
\gcd(a_0, a_1) & \mid d - a_2\\
\gcd(a_0, a_1) & \mid d - a_3.
\end{align*}
Therefore $\gcd(a_0, a_1) \mid a_2, a_3$, and condition~(ii) gives $\gcd(a_0, a_1) = 1$.
So, (v) is satisfied.
\end{proof}
\end{myproposition}

\begin{mytheorem}{\ref*{theorem:series-check-one}} \label{proof:series-check-one}
\StateTheoremSeriesCheckOne

\begin{proof}
By Proposition~\ref{proposition:four-conds-sufficient},
it is sufficient to show the tuple $(a_0, a_1, a_2, a_3, d)$ satisfies conditions (i)--(iv).
Note that (iii) clearly holds for all the classes (1)--(6).
By Proposition~\ref{proposition:conds-imply-class}, it suffices to consider tuples of class 1--6.
\smallskip

\emph{Class 1}: $I = a_0 + a_1$, giving
\[
d = a_2 + a_3.
\]
From Definition~\ref{definition:corresponding-series}, there exist $x$ and $y$ such that
\begin{align*}
m & = \lcm(a_0, a_1)\\
(a_0, a_1, a_2, a_3, d) & = (b_0, b_1, b_2+xm, b_3+ym, e+(x+y)m)
\end{align*}
and the tuple is ordered, that is, $a_0 \leq a_1 \leq a_2 \leq a_3 \leq d$.

First, we check condition~(i). Since $i < j$, we find $i \in \{0,1,2\}$.
\begin{itemize}
\item If $i \in \{0,1\}$, then $a_i = b_i$ and $a_i \mid m$, giving $\gcd(a_i, a_j) = \gcd(b_i, b_j)$.
By assumption, the tuple $(b_0, b_1, b_2, b_3, e)$ satisfies conditions (i)--(iv), so $\gcd(b_i, b_j) \mid e$, giving $\gcd(a_i, a_j) \mid d$.
\item For $(i,j) = (2,3)$, we note $d = a_2 + a_3$.
\end{itemize}

Next, we check condition~(ii).
Since $i < j < k$, we find $i \in \{0,1\}$, giving $a_i = b_i$ and $a_i \mid m$.
Therefore, $\gcd(a_i, a_j, a_k) = \gcd(b_i, b_j, b_k)$, and by assumption $\gcd(b_i, b_j, b_k) = 1$.

Finally, we check condition~(iv).
\begin{itemize}
\item If $i \in \{0,1\}$, then $a_i = b_i$ and $a_i \mid m$.
By assumption, there exists $j$ such that $b_i \mid e - b_j$, giving $a_i \mid d - a_j$.
\item If $i \in \{2,3\}$, we note $d = a_2 + a_3$.
\end{itemize}
\smallskip

\emph{Classes 2--3}: $I = a_t + a_2$, where $t \in \{0,1\}$.
We can define $u$ such that
\begin{align*}
\{t,u\} & = \{0,1\}\\
d & = a_u + a_3.
\end{align*}
From Definition~\ref{definition:corresponding-series}, there exist $x$ and $y$ such that
\begin{align*}
m & = \lcm(a_0, a_1, a_2)\\
(a_0, a_1, a_2, a_3, d) & = (b_0, b_1, b_2, b_3+xm, e+xm).
\end{align*}

First, we check condition~(i).
Since $i < j$, we find $i \in \{0,1,2\}$, giving $a_i = b_i$ and $a_i \mid m$.
We have $\gcd(a_i, a_j) = \gcd(b_i, b_j)$.
By assumption $\gcd(b_i, b_j) \mid e$, giving $\gcd(a_i, a_j) \mid d$.

For condition~(ii), we similarly find $\gcd(a_i, a_j, a_k) = \gcd(b_i, b_j, b_k) = 1$.

Finally, we consider condition~(iv).
It holds for $i = 3$, since $d = a_u + a_3$.
If $i \neq 3$, then by assumption there exists $j$ such that $b_i \mid e - b_j$.
Since $a_i = b_i$ and $a_i \mid m$, we find $a_i \mid d - a_j$.
\smallskip

\emph{Classes 4--5}: $I = a_t + \frac{a_u}{2}$, where $t$ and $u$ are such that
\[
\{t,u\} = \{0,1\}.
\]
From Definition~\ref{definition:corresponding-series}, there exist $x$ and $y$ such that
\begin{align*}
a_3 & = \frac{a_u}{2} + a_2\\
d & = 2 a_3 = a_u + 2 a_2\\
m & = \lcm\left(a_t, \frac{a_u}{2}\right)\\
(a_0, a_1, a_2, a_3, d) & = (b_0, b_1, b_2+xm, b_3+xm, e+2xm).
\end{align*}

First, we check condition~(i). We have $i < j$.
\begin{itemize}
\item If $j = 3$, we note $d = 2 a_3$.
\item If $(i,j) = (u,2)$, we note $d = a_u + 2 a_2$.
\item If $i = t$ or $j = t$, then $a_t = b_t$ and $a_t \mid m$, giving $\gcd(a_i, a_j) = \gcd(b_i, b_j)$.
By assumption $\gcd(b_i, b_j) \mid e$, giving $\gcd(a_i, a_j) \mid d$.
\end{itemize}

Next, we check condition~(ii).
\begin{itemize}
\item If $i = t$ or $j = t$, we find $a_t = b_t$ and $a_t \mid m$, giving $\gcd(a_i, a_j, a_k) = \gcd(b_i, b_j, b_k) = 1$.
\item If $(i,j,k) = (u,2,3)$, then using $a_3 = \frac{a_u}{2} + a_2$, we find $\gcd(a_u, a_2, a_3) = \gcd(\frac{a_u}2, a_2)$.
We have $a_u = b_u$ and $\frac{a_u}2 \mid m$, therefore $\gcd(\frac{a_u}2, a_2) = \gcd(\frac{b_u}2, b_2) = \gcd(b_u, b_2, b_3)$.
By assumption $\gcd(b_u, b_2, b_3) = 1$, giving $\gcd(a_u, a_2, a_3) = 1$.
\end{itemize}

Finally, we check condition~(iv).
\begin{itemize}
\item If $i = 3$, we note $d = 2 a_3$.
\item If $i = 2$, we note $d = a_u + 2 a_2$.
\item If $i = t$, then $a_i = b_i$ and $a_i \mid m$.
By assumption there exists $j$ such that $b_i \mid e - b_j$, therefore $a_i \mid d - a_j$.
\item If $i = u$, then
\begin{align*}
a_u & = b_u\\
a_u & \mid 2m.
\end{align*}
From Definition~\ref{definition:corresponding-series}
\[
e = 2 b_3 = b_u + 2 b_2.
\]
By assumption there exists $j$ such that $b_u \mid e - b_j$.
We show there exists $k$ such that $a_u \mid e - a_k$.
\begin{itemize}
\item If $j \in \{0,1\}$, then $a_j = b_j$, giving $a_u \mid e - a_j$.
\item If $j \in \{2,3\}$, then either $b_u \mid b_2$ or $b_u \mid b_3$, giving $b_u \mid e$.
This implies $a_u \mid e - a_u$.
\end{itemize}
Now, since $a_u \mid 2m$, we find $a_u \mid d - a_k$.
So, condition~(iv) is satisfied.
\end{itemize}
\smallskip

\emph{Class 6}: there exist $a$ and $k$ such that
\begin{align*}
1 & \leq k < I\\
a & \geq I + k\\
(a_0, a_1, a_2, a_3, d) & = (I-k, I+k, a, a+k, 2a+I+k).
\end{align*}
We have
\[
d = a_0 + 2 a_3 = a_1 + 2 a_2.
\]
From Definition~\ref{definition:corresponding-series}, we have
\begin{align*}
m & = \lcm(a_0, a_1, k)\\
(a_0, a_1, a_2, a_3, d) & = (b_0, b_1, b_2+xm, b_3+xm, e+2xm).
\end{align*}

First, we check condition~(i). Since $i < j$, we have $i \in \{0,1,2\}$.
\begin{itemize}
\item If $i \in \{0,1\}$, then $a_i = b_i$ and $a_i \mid m$, giving $\gcd(a_i, a_j) = \gcd(b_i, b_j)$.
By assumption $\gcd(b_i, b_j) \mid e$, giving $\gcd(a_i, a_j) \mid d$.
\item If $(i,j) = (2,3)$, then $\gcd(a_2, a_3) = \gcd(a_2, k)$.
Since $k \mid m$, we find $\gcd(a_2, k) = \gcd(b_2, k) = \gcd(b_2, b_3)$.
By assumption $\gcd(b_2, b_3) \mid e$, giving $\gcd(a_2, a_3) \mid d$.
\end{itemize}

Next, we check condition~(ii).
Since $i < j < k$, we have $i \in \{0,1\}$, giving $a_i = b_i$ and $a_i \mid m$.
We find $\gcd(a_i, a_j, a_k) = \gcd(b_i, b_j, b_k) = 1$.

Finally, we check condition~(iv).
\begin{itemize}
\item If $i = 3$, we note $d = a_0 + 2 a_3$.
\item If $i = 2$, we note $d = a_1 + 2 a_2$.
\item If $i \in \{0,1\}$, then $a_i = b_i$ and $a_i \mid m$. By assumption $b_i \mid e - b_j$, giving $a_i \mid d - a_j$.
\end{itemize}
\end{proof}
\end{mytheorem}

\bigskip
\section{Small Index Cases} \label{section:small-index-cases}

In this section, we give the complete lists of quasi-smooth well-formed hypersurfaces for indices $I = 1,2,\ldots,6$.
The parameters $x$ and $y$ are non-negative integers with $x \leq y$.
We first list the two-parameter series, then one-parameter series and lastly sporadic cases.

The tables for indices $1$ and~$2$ were known, from \cite[Theorem~8]{JoKo01b} and \cite[Corollary~1.13]{ChSh-Zoo} respectively.
The tables for indices $3$ to $6$ are new.
The author has computed the lists for all $I \leq 100$.
The lists grow as the cube of the index and the computation time grows as the fifth power of~$I$.

Differences from~\cite[Corollary~1.13]{ChSh-Zoo}: There is a misprint in the list for $I = 2$ in~\cite{ChSh-Zoo},
namely, the second occurrence of $(3, 4, 6, 7, 18)$ should instead be $(3, 4, 5, 7, 17)$.
In the list below, the tuple $(3, 4, 5, 7, 17)$ is contained in the series $(3, 3x + 4, 3x + 5, 6x + 7, 12x + 17)$
which has been extended to include $x = 0$.
The tuple $(1, 1, 2, 2, 4)$ is contained in the series $(1, 1, x+1, y+1, x+y+2)$.


\begin{longtable}{|c|c||c|c||c|c|}
\caption{Index 1} \label{table:ind-1}\\
\hline
$(a_0,a_1,a_2,a_3)$ & $d$ & $(a_0,a_1,a_2,a_3)$ & $d$ & $(a_0,a_1,a_2,a_3)$ & $d$
\endhead
\hline

$(2,2x+3,2x+3,4x+5)$ & $8x+12$ &
$(1,1,1,1)$ & $3$ &
$(1,1,1,2)$ & $4$\\
\hline

$(1,1,2,3)$ & $6$ &
$(1,2,3,5)$ & $10$ &
$(1,3,5,7)$ & $15$\\
\hline

$(1,3,5,8)$ & $16$ &
$(2,3,5,9)$ & $18$ &
$(3,3,5,5)$ & $15$\\
\hline

$(3,5,7,11)$ & $25$ &
$(3,5,7,14)$ & $28$ &
$(3,5,11,18)$ & $36$\\
\hline

$(5,14,17,21)$ & $56$ &
$(5,19,27,31)$ & $81$ &
$(5,19,27,50)$ & $100$\\
\hline

$(7,11,27,37)$ & $81$ &
$(7,11,27,44)$ & $88$ &
$(9,15,17,20)$ & $60$\\
\hline

$(9,15,23,23)$ & $69$ &
$(11,29,39,49)$ & $127$ &
$(11,49,69,128)$ & $256$\\
\hline

$(13,23,35,57)$ & $127$ &
$(13,35,81,128)$ & $256$ & & \\
\hline
\end{longtable}

\begin{longtable}{|c|c||c|c|}
\caption{Index 2} \label{table:ind-2}\\
\hline
$(a_0,a_1,a_2,a_3)$ & $d$ & $(a_0,a_1,a_2,a_3)$ & $d$
\endhead
\hline

$(1,1,x+1,y+1)$ & $x+y+2$ &
$(1,2,x+2,x+3)$ & $2x+6$\\
\hline

$(1,3,3x+3,3x+4)$ & $6x+9$ &
$(1,3,3x+4,3x+5)$ & $6x+11$\\
\hline

$(3,3x+3,3x+4,3x+4)$ & $9x+12$ &
$(3,3x+4,3x+5,3x+5)$ & $9x+15$\\
\hline

$(3,3x+4,3x+5,6x+7)$ & $12x+17$ &
$(3,3x+4,6x+7,9x+9)$ & $18x+21$\\
\hline

$(3,3x+4,6x+7,9x+12)$ & $18x+24$ &
$(4,2x+5,2x+5,4x+8)$ & $8x+20$\\
\hline

$(4,2x+5,4x+10,6x+13)$ & $12x+30$ &
$(1,3,4,6)$ & $12$\\
\hline

$(1,4,5,7)$ & $15$ &
$(1,4,5,8)$ & $16$\\
\hline

$(1,4,6,9)$ & $18$ &
$(1,5,7,11)$ & $22$\\
\hline

$(1,6,9,13)$ & $27$ &
$(1,6,10,15)$ & $30$\\
\hline

$(1,7,12,18)$ & $36$ &
$(1,8,13,20)$ & $40$\\
\hline

$(1,9,15,22)$ & $45$ &
$(2,3,4,5)$ & $12$\\
\hline

$(2,3,4,7)$ & $14$ &
$(3,4,5,10)$ & $20$\\
\hline

$(3,4,6,7)$ & $18$ &
$(3,4,10,15)$ & $30$\\
\hline

$(5,13,19,22)$ & $57$ &
$(5,13,19,35)$ & $70$\\
\hline

$(6,9,10,13)$ & $36$ &
$(7,8,19,25)$ & $57$\\
\hline

$(7,8,19,32)$ & $64$ &
$(9,12,13,16)$ & $48$\\
\hline

$(9,12,19,19)$ & $57$ &
$(9,19,24,31)$ & $81$\\
\hline

$(10,19,35,43)$ & $105$ &
$(11,21,28,47)$ & $105$\\
\hline

$(11,25,32,41)$ & $107$ &
$(11,25,34,43)$ & $111$\\
\hline

$(11,43,61,113)$ & $226$ &
$(13,18,45,61)$ & $135$\\
\hline

$(13,20,29,47)$ & $107$ &
$(13,20,31,49)$ & $111$\\
\hline

$(13,31,71,113)$ & $226$ &
$(14,17,29,41)$ & $99$\\
\hline
\end{longtable}

\begin{longtable}{|c|c||c|c|}
\caption{Index 3} \label{table:ind-3}\\
\hline
$(a_0,a_1,a_2,a_3)$ & $d$ & $(a_0,a_1,a_2,a_3)$ & $d$
\endhead
\hline

$(1,2,2x+3,2y+3)$ & $2(x+y)+6$ &
$(1,1,2,2x+3)$ & $2x+4$\\
\hline

$(1,5,10x+5,10x+7)$ & $20x+15$ &
$(1,5,10x+7,10x+9)$ & $20x+19$\\
\hline

$(1,7,9,13)$ & $27$ &
$(1,7,9,14)$ & $28$\\
\hline

$(1,9,13,20)$ & $40$ &
$(1,13,22,33)$ & $66$\\
\hline

$(1,14,23,35)$ & $70$ &
$(1,15,25,37)$ & $75$\\
\hline

$(5,7,11,13)$ & $33$ &
$(5,7,11,20)$ & $40$\\
\hline

$(11,21,29,37)$ & $95$ &
$(11,37,53,98)$ & $196$\\
\hline

$(13,17,27,41)$ & $95$ &
$(13,27,61,98)$ & $196$\\
\hline

$(15,19,43,74)$ & $148$& & \\
\hline
\end{longtable}

\begin{longtable}{|c|c|}
\caption{Index 4} \label{table:ind-4}\\
\hline
$(a_0,a_1,a_2,a_3)$ & $d$
\endhead
\hline

$(1,3,3x+4,3y+5)$ & $3(x+y)+9$\\
\hline

$(1,3,3x+5,3y+5)$ & $3(x+y)+10$\\
\hline

$(1,3,3x+5,3y+7)$ & $3(x+y)+12$\\
\hline

$(2,2,2x+3,2y+3)$ & $2(x+y)+6$\\
\hline

$(1,1,3,3x+5)$ & $3x+6$\\
\hline

$(1,2,2,2x+3)$ & $2x+4$\\
\hline

$(1,2,3,6x+4)$ & $6x+6$\\
\hline

$(1,2,3,6x+5)$ & $6x+7$\\
\hline

$(1,2,3,6x+7)$ & $6x+9$\\
\hline

$(1,2,3,6x+8)$ & $6x+10$\\
\hline

$(1,7,21x+7,21x+10)$ & $42x+21$\\
\hline

$(1,7,21x+10,21x+13)$ & $42x+27$\\
\hline

$(1,7,21x+14,21x+17)$ & $42x+35$\\
\hline

$(1,7,21x+17,21x+20)$ & $42x+41$\\
\hline

$(2,3,3x+4,3x+5)$ & $6x+10$\\
\hline

$(2,3,3x+5,3x+6)$ & $6x+12$\\
\hline

$(2,4,2x+5,2x+7)$ & $4x+14$\\
\hline

$(2,6,6x+9,6x+11)$ & $12x+24$\\
\hline

$(3,5,15x+5,15x+6)$ & $30x+15$\\
\hline

$(3,5,15x+10,15x+11)$ & $30x+25$\\
\hline

$(3,5,15x+11,15x+12)$ & $30x+27$\\
\hline

$(3,5,15x+16,15x+17)$ & $30x+37$\\
\hline

$(6,6x+9,6x+11,6x+11)$ & $18x+33$\\
\hline

$(6,6x+11,12x+20,18x+27)$ & $36x+60$\\
\hline

$(6,6x+11,12x+20,18x+33)$ & $36x+66$\\
\hline

$(1,10,13,19)$ & $39$\\
\hline

$(1,10,13,20)$ & $40$\\
\hline

$(1,13,19,29)$ & $58$\\
\hline

$(1,14,21,31)$ & $63$\\
\hline

$(1,19,32,48)$ & $96$\\
\hline

$(1,20,33,50)$ & $100$\\
\hline

$(1,21,35,52)$ & $105$\\
\hline

$(2,7,10,15)$ & $30$\\
\hline

$(2,9,12,17)$ & $36$\\
\hline

$(5,6,8,9)$ & $24$\\
\hline

$(5,6,8,15)$ & $30$\\
\hline

$(9,11,12,17)$ & $45$\\
\hline

$(10,13,25,31)$ & $75$\\
\hline

$(11,17,20,27)$ & $71$\\
\hline

$(11,17,24,31)$ & $79$\\
\hline

$(11,31,45,83)$ & $166$\\
\hline

$(13,14,19,29)$ & $71$\\
\hline

$(13,14,23,33)$ & $79$\\
\hline

$(13,23,51,83)$ & $166$\\
\hline
\end{longtable}

\begin{longtable}{|c|c||c|c|}
\caption{Index 5} \label{table:ind-5}\\
\hline
$(a_0,a_1,a_2,a_3)$ & $d$ & $(a_0,a_1,a_2,a_3)$ & $d$
\endhead
\hline

$(1,4,4x+5,4y+7)$ & $4(x+y)+12$ &
$(1,4,4x+7,4y+9)$ & $4(x+y)+16$\\
\hline

$(2,3,6x+5,6y+7)$ & $6(x+y)+12$ &
$(2,3,6x+7,6y+7)$ & $6(x+y)+14$\\
\hline

$(2,3,6x+7,6y+11)$ & $6(x+y)+18$ &
$(1,1,4,4x+7)$ & $4x+8$\\
\hline

$(1,2,3,6x+5)$ & $6x+6$ &
$(1,2,3,6x+7)$ & $6x+8$\\
\hline

$(1,3,4,12x+5)$ & $12x+8$ &
$(1,3,4,12x+9)$ & $12x+12$\\
\hline

$(1,3,4,12x+13)$ & $12x+16$ &
$(1,9,36x+9,36x+13)$ & $72x+27$\\
\hline

$(1,9,36x+13,36x+17)$ & $72x+35$ &
$(1,9,36x+27,36x+31)$ & $72x+63$\\
\hline

$(1,9,36x+31,36x+35)$ & $72x+71$ &
$(3,7,42x+7,42x+9)$ & $84x+21$\\
\hline

$(3,7,42x+23,42x+25)$ & $84x+53$ &
$(3,7,42x+35,42x+37)$ & $84x+77$\\
\hline

$(3,7,42x+37,42x+39)$ & $84x+81$ &
$(1,13,17,25)$ & $51$\\
\hline

$(1,13,17,26)$ & $52$ &
$(1,17,25,38)$ & $76$\\
\hline

$(1,25,42,63)$ & $126$ &
$(1,26,43,65)$ & $130$\\
\hline

$(1,27,45,67)$ & $135$ &
$(6,7,9,10)$ & $27$\\
\hline

$(11,13,19,25)$ & $63$ &
$(11,25,37,68)$ & $136$\\
\hline

$(13,19,41,68)$ & $136$ & & \\
\hline
\end{longtable}

\begin{longtable}{|c|c||c|c|}
\caption{Index 6} \label{table:ind-6}\\
\hline
$(a_0,a_1,a_2,a_3)$ & $d$ & $(a_0,a_1,a_2,a_3)$ & $d$
\endhead
\hline

$(1,5,5x+6,5y+9)$ & $5(x+y)+15$ &
$(1,5,5x+7,5y+8)$ & $5(x+y)+15$\\
\hline

$(1,5,5x+7,5y+9)$ & $5(x+y)+16$ &
$(1,5,5x+8,5y+8)$ & $5(x+y)+16$\\
\hline

$(1,5,5x+8,5y+12)$ & $5(x+y)+20$ &
$(1,5,5x+9,5y+11)$ & $5(x+y)+20$\\
\hline

$(1,5,5x+9,5y+12)$ & $5(x+y)+21$ &
$(2,4,4x+5,4y+5)$ & $4(x+y)+10$\\
\hline

$(2,4,4x+5,4y+7)$ & $4(x+y)+12$ &
$(2,4,4x+7,4y+7)$ & $4(x+y)+14$\\
\hline

$(2,4,4x+7,4y+9)$ & $4(x+y)+16$ &
$(3,3,3x+4,3y+5)$ & $3(x+y)+9$\\
\hline

$(3,3,3x+5,3y+7)$ & $3(x+y)+12$ &
$(1,1,5,5x+9)$ & $5x+10$\\
\hline

$(1,2,4,4x+5)$ & $4x+6$ &
$(1,2,4,4x+7)$ & $4x+8$\\
\hline

$(1,2,5,10x+8)$ & $10x+10$ &
$(1,2,5,10x+9)$ & $10x+11$\\
\hline

$(1,2,5,10x+13)$ & $10x+15$ &
$(1,2,5,10x+14)$ & $10x+16$\\
\hline

$(1,3,3,3x+5)$ & $3x+6$ &
$(1,3,5,15x+7)$ & $15x+10$\\
\hline

$(1,3,5,15x+8)$ & $15x+11$ &
$(1,3,5,15x+12)$ & $15x+15$\\
\hline

$(1,3,5,15x+13)$ & $15x+16$ &
$(1,3,5,15x+17)$ & $15x+20$\\
\hline

$(1,3,5,15x+18)$ & $15x+21$ &
$(1,4,5,20x+6)$ & $20x+10$\\
\hline

$(1,4,5,20x+7)$ & $20x+11$ &
$(1,4,5,20x+11)$ & $20x+15$\\
\hline

$(1,4,5,20x+12)$ & $20x+16$ &
$(1,4,5,20x+16)$ & $20x+20$\\
\hline

$(1,4,5,20x+17)$ & $20x+21$ &
$(1,4,5,20x+21)$ & $20x+25$\\
\hline

$(1,4,5,20x+22)$ & $20x+26$ &
$(1,11,55x+11,55x+16)$ & $110x+33$\\
\hline

$(1,11,55x+16,55x+21)$ & $110x+43$ &
$(1,11,55x+22,55x+27)$ & $110x+55$\\
\hline

$(1,11,55x+27,55x+32)$ & $110x+65$ &
$(1,11,55x+33,55x+38)$ & $110x+77$\\
\hline

$(1,11,55x+38,55x+43)$ & $110x+87$ &
$(1,11,55x+44,55x+49)$ & $110x+99$\\
\hline

$(1,11,55x+49,55x+54)$ & $110x+109$ &
$(2,3,3,6x+4)$ & $6x+6$\\
\hline

$(2,3,3,6x+7)$ & $6x+9$ &
$(2,3,4,12x+5)$ & $12x+8$\\
\hline

$(2,3,4,12x+7)$ & $12x+10$ &
$(2,3,4,12x+9)$ & $12x+12$\\
\hline

$(2,3,4,12x+11)$ & $12x+14$ &
$(2,3,4,12x+13)$ & $12x+16$\\
\hline

$(2,3,4,12x+15)$ & $12x+18$ &
$(2,5,5x+8,5x+9)$ & $10x+18$\\
\hline

$(2,5,5x+9,5x+10)$ & $10x+20$ &
$(2,8,4x+9,4x+13)$ & $8x+26$\\
\hline

$(2,10,20x+15,20x+19)$ & $40x+40$ &
$(2,10,20x+25,20x+29)$ & $40x+60$\\
\hline

$(5,7,35x+8,35x+9)$ & $70x+23$ &
$(5,7,35x+14,35x+15)$ & $70x+35$\\
\hline

$(5,7,35x+28,35x+29)$ & $70x+63$ &
$(5,7,35x+29,35x+30)$ & $70x+65$\\
\hline

$(8,4x+9,4x+11,4x+13)$ & $12x+35$ &
$(9,3x+11,3x+14,6x+19)$ & $12x+47$\\
\hline

$(1,16,21,31)$ & $63$ &
$(1,16,21,32)$ & $64$\\
\hline

$(1,21,31,47)$ & $94$ &
$(1,22,33,49)$ & $99$\\
\hline

$(1,31,52,78)$ & $156$ &
$(1,32,53,80)$ & $160$\\
\hline

$(1,33,55,82)$ & $165$ &
$(2,13,18,27)$ & $54$\\
\hline

$(2,15,20,29)$ & $60$ &
$(3,7,8,12)$ & $24$\\
\hline

$(7,10,15,19)$ & $45$ &
$(11,19,29,53)$ & $106$\\
\hline

$(13,15,31,53)$ & $106$ & & \\
\hline
\end{longtable}

\appendix

\bigskip
\section{Tables} \label{section:infinite-series}

Tables \ref{table:inf} and~\ref{table:spor} are taken from \cite[Appendix B]{ChSh-Zoo}.
They contain one-parameter infinite series and sporadic cases respectively of~values of $(a_0,a_1,a_2,a_3,d,I)$.
The~last columns represent the~cases in~\cite{YauYu03} from which the~sextuples $(a_0,a_1,a_2,a_3,d,I)$
originate\footnote{Note that sometimes a sextuple $(a_0,a_1,a_2,a_3,d,I)$ originates from several cases in~\cite{YauYu03}.}.
The parameter $n$ is any positive integer.

Differences from~\cite[Appendix B]{ChSh-Zoo}:
the tuple $(3, 3, 4, 4, 12)$ with $I = 2$ has been removed from Table~\ref{table:spor},
as it already appears in the series $(3, 3n, 3n + 1, 3n + 1, 9n + 3)$ with $I = 2$ in Table~\ref{table:inf}.

\begin{longtable}{|c|c|c|c|}
\caption{Infinite series} \label{table:inf}\\ \hline
$(a_0,a_1,a_2,a_3)$ & $d$ & $I$ & Source
\endhead
\hline

$(1, 3n-2, 4n-3, 6n-5)$ & $12n-9$ & $n$ & VII.2(3) \\
\hline

$(1, 3n-2, 4n-3, 6n-4)$ & $12n-8$ & $n$ & II.2(2) \\
\hline

$(1, 4n-3, 6n-5, 9n-7)$ & $18n-14$ & $n$ &VII.3(1) \\
\hline

$(1, 6n-5, 10n-8, 15n-12)$ & $30n-24$ & $n$ &  III.1(4) \\
\hline

$(1, 6n-4, 10n-7, 15n-10)$ & $30n-20$ & $n$ & III.2(2) \\
\hline

$(1, 6n-3, 10n-5, 15n-8)$ & $30n-15$ & $n$ & III.2(4) \\
\hline

$(1, 8n-2, 12n-3, 18n-5)$ & $36n-9$ & $2n$ & IV.3(3) \\
\hline

$(2, 6n-3, 8n-4, 12n-7)$ & $24n-12$ & $2n$ & II.2(4)\\
\hline

$(2, 6n+1, 8n+2, 12n+3)$ & $24n+6$ & $2n+2$ & II.2(1)\\
\hline

$(3, 6n+1, 6n+2, 9n+3)$ & $18n+6$ & $3n+3$ & II.2(1)\\
\hline

$(7, 28n-18, 42n-27, 63n-44)$ & $126n-81$ & $7n-1$ & XI.3(14)\\
\hline

$(7, 28n-17, 42n-29, 63n-40)$ & $126n-80$ & $7n+1$ & X.3(1)\\
\hline

$(7, 28n-13, 42n-23, 63n-31)$ & $126n-62$ & $7n+2$ & X.3(1)\\
\hline

$(7, 28n-10, 42n-15, 63n-26)$ & $126n-45$ & $7n+1$ & XI.3(14)\\
\hline

$(7, 28n-9, 42n-17, 63n-22)$ & $126n-44$ & $7n+3$ & X.3(1)\\
\hline

$(7, 28n-6, 42n-9, 63n-17)$ & $126n-27$ & $7n+2$ & XI.3(14)\\
\hline

$(7, 28n-5, 42n-11, 63n-13)$ & $126n-26$ & $7n+4$ & X.3(1)\\
\hline

$(7, 28n-2, 42n-3, 63n-8)$ & $126n-9$ & $7n+3$ & XI.3(14)\\
\hline

$(7, 28n-1, 42n-5, 63n-4)$ & $126n-8$ & $7n+5$ & X.3(1)\\
\hline

$(7, 28n+2, 42n+3, 63n+1)$ & $126n+9$ & $7n+4$ & XI.3(14)\\
\hline

$(7, 28n+3, 42n+1, 63n+5)$ & $126n+10$ & $7n+6$ & X.3(1)\\
\hline

$(7, 28n+6, 42n+9, 63n+10)$ & $126n+27$ & $7n+5$ & XI.3(14)\\
\hline

$(2, 2n+1, 2n+1, 4n+1)$ & $8n+4$ & $1$ & II.3(4) \\
\hline

$(3, 3n, 3n+1, 3n+1)$ & $9n+3$ & $2$ & III.5(1) \\
\hline

$(3, 3n+1, 3n+2, 3n+2)$ & $9n+6$ & $2$ & II.5(1) \\
\hline

$(3, 3n+1, 3n+2, 6n+1)$ & $12n+5$ & $2$ & XVIII.2(2)\\
\hline

$(3, 3n+1, 6n+1, 9n)$ & $18n+3$ & $2$ & VII.3(2)\\
\hline

$(3, 3n+1, 6n+1, 9n+3)$ & $18n+6$ & $2$ & II.2(2) \\
\hline

$(4, 2n+3, 2n+3, 4n+4)$ & $8n+12$ & $2$ &  V.3(4) \\
\hline

$(4, 2n+3, 4n+6, 6n+7)$ & $12n+18$ & $2$ & XII.3(17) \\
\hline

$(6, 6n+3, 6n+5, 6n+5)$ & $18n+15$ & $4$ & III.5(1) \\
\hline

$(6, 6n+5, 12n+8, 18n+9)$ & $36n+24$ & $4$ & VII.3(2)\\
\hline

$(6, 6n+5, 12n+8, 18n+15)$ & $36n+30$ & $4$ & IV.3(1) \\
\hline

$(8, 4n+5, 4n+7, 4n+9)$ & $12n+23$ & $6$ & XIX.2(2)\\
\hline

$(9, 3n+8, 3n+11, 6n+13)$ & $12n+35$ & $6$ & XIX.2(2)\\
\hline
\end{longtable}

\begin{longtable}{|c|c|c|c||c|c|c|c|}
\caption{Sporadic cases} \label{table:spor}\\
\hline
$(a_0,a_1,a_2,a_3)$ & $d$ & $I$ & Source & $(a_0,a_1,a_2,a_3)$ & $d$ & $I$ & Source
\endhead
\hline

$(1, 3, 5, 8)$ & $16$ & $1$ & VIII.3(5) &
$(2, 3, 5, 9)$ & $18$ & $1$ & II.2(3)\\
\hline

$(3, 3, 5, 5)$ & $15$ & $1$ & I.19 &
$(3, 5, 7, 11)$ & $25$ & $1$ & X.2(3)\\
\hline

$(3, 5, 7, 14)$ & $28$ & $1$ & VII.4(4)&
$(3, 5, 11, 18)$ & $36$ & $1$ & VII.3(1)\\
\hline

$(5, 14, 17, 21)$ & $56$ & $1$ & XI.3(8)&
$(5, 19, 27, 31)$ & $81$ & $1$ & X.3(3)\\
\hline

$(5, 19, 27, 50)$ & $100$ & $1$ & VII.3(3)&
$(7, 11, 27, 37)$ & $81$ & $1$ & X.3(4)\\
\hline

$(7, 11, 27, 44)$ & $88$ & $1$ & VII.3(5)&
$(9, 15, 17, 20)$ & $60$ & $ 1$ & VII.6(3)\\
\hline

$(9, 15, 23, 23)$ & $69$ & $1$ & III.5(1) &
$(11, 29, 39, 49)$ & $127$ & $1$ & XIX.2(2)\\
\hline

$(11, 49, 69, 128)$ & $256$ & $1$ & X.3(1)&
$(13, 23, 35, 57)$ & $127$ & $1$ & XIX.2(2)\\
\hline

$(13, 35, 81, 128)$ & $256$ & $1$ & X.3(2)&
$(1, 3, 4, 6)$ & $12$ & $2$ & I.3 \\
\hline

$(1, 4, 6, 9)$ & $18$ & $2$ & IV.3(3) &
$(1, 6, 10, 15)$ & $30$ & $2$ & I.4\\
\hline

$(2, 3, 4, 7)$ & $14$ & $2$ & IX.3(1)&
$(3, 4, 5, 10)$ & $20$ & $2$ & II.3(2)\\
\hline

$(3, 4, 6, 7)$ & $18$ & $2$ & VII.3(10) &
$(3, 4, 10, 15)$ & $30$ & $2$ & II.2(3)\\
\hline

$(5, 13, 19, 22)$ & $57$ & $2$ & X.3(3) &
$(5, 13, 19, 35)$ & $70$ & $2$ & VII.3(3)\\
\hline

$(6, 9, 10, 13)$ & $36$ & $2$ & VII.3(8) &
$(7, 8, 19, 25)$ & $57$ & $2$ & X.3(4)\\
\hline

$(7, 8, 19, 32)$ & $64$ & $2$ & VII.3(3) &
$(9, 12, 13, 16)$ & $48$ & $2$ & VII.6(2)\\
\hline

$(9, 12, 19, 19)$ & $57$ & $2$ & III.5(1) &
$(9, 19, 24, 31)$ & $81$ & $2$ & XI.3(20)\\
\hline

$(10, 19, 35, 43)$ & $105$ & $2$ & XI.3(18) &
$(11, 21, 28, 47)$ & $105$ & $2$ & XI.3(16)\\
\hline

$(11, 25, 32, 41)$ & $107$ & $2$ & XIX.3(1) &
$(11, 25, 34, 43)$ & $111$ & $2$ & XIX.2(2)\\
\hline

$(11, 43, 61, 113)$ & $226$ & $2$ & X.3(1) &
$(13, 18, 45, 61)$ & $135$ & $2$ & XI.3(14)\\
\hline

$(13, 20, 29, 47)$ & $107$ & $2$ & XIX.3(1) &
$(13, 20, 31, 49)$ & $111$ & $2$ & XIX.2(2)\\
\hline

$(13, 31, 71, 113)$ & $226$ & $2$ & X.3(2) &
$(14, 17, 29, 41)$ & $99$ & $2$ & XIX.2(3)\\
\hline

$(5, 7, 11, 13)$ & $33$ & $3$ & X.3(3) &
$(5, 7, 11, 20)$ & $40$ & $3$ & VII.3(3)\\
\hline

$(11, 21, 29, 37)$ & $95$ & $3$ & XIX.2(2) &
$(11, 37, 53, 98)$ & $196$ & $3$ & X.3(1)\\
\hline

$(13, 17, 27, 41)$ & $95$ & $3$ & XIX.2(2) &
$(13, 27, 61, 98)$ & $196$ & $3$ & X.3(2)\\
\hline

$(15, 19, 43, 74)$ & $148$ & $3$ & X.3(1) &
$(5, 6, 8, 9)$ & $24$ & $4$ & VII.3(2)\\
\hline

$(5, 6, 8, 15)$ & $30$ & $4$ & IV.3(1) &
$(9, 11, 12, 17)$ & $45$ & $4$ & XI.3(20)\\
\hline

$(10, 13, 25, 31)$ & $75$ & $4$ & XI.3(14) &
$(11, 17, 20, 27)$ & $71$ & $4$ & XIX.3(1)\\
\hline

$(11, 17, 24, 31)$ & $79$ & $4$ & XIX.2(2) &
$(11, 31, 45, 83)$ & $166$ & $4$ & X.3(1)\\
\hline

$(13, 14, 19, 29)$ & $71$ & $4$ & XIX.3(1) &
$(13, 14, 23, 33)$ & $79$ & $4$ & XIX.2(2)\\
\hline

$(13, 23, 51, 83)$ & $166$ & $4$ & X.3(2) &
$(6, 7, 9, 10)$ & $27$ & $5$ & XI.3(14)\\
\hline

$(11, 13, 19, 25)$ & $63$ & $5$ & XIX.2(2) &
$(11, 25, 37, 68)$ & $136$ & $5$ & X.3(1)\\
\hline

$(13, 19, 41, 68)$ & $136$ & $5$ & X.3(2) &
$(11, 19, 29, 53)$ & $106$ & $6$ & X.3(1)\\
\hline

$(13, 15, 31, 53)$ & $106$ & $6$ & X.3(2) &
$(11, 13, 21, 38)$ & $76$ & $7$ & X.3(1)\\
\hline
\end{longtable}

\bigskip
\section{Source Code} \label{section:source-code}
The computer code below classifies the hypersurfaces of index $I$.
For simplicity, the tuples from the tables are left out.
The full program and source code are available from the author.
It is written in the functional programming language Haskell.

\begin{lstlisting}
-- Classifying quasi-smooth well-formed weighted hypersurfaces.
-- Erik Paemurru

data Tuple = Quint Int Int Int Int Int Int Int deriving( Eq,Ord,Show )
--                (a0, a1, a2, a3, dd, mm, cc)
-- dd - degree
-- mm - series modulo-number (lcm of smaller weights)
-- cc - series class-number (from 1 to 6)

main = do
  putStr ("Enter index, for which to solve:\n" ++ "Index = ")
  strL <- getLine
  mapM_ putStrLn (map show (solve (read strL :: Int)))

-- The 'solve' function classifies the hypersurfaces. It selects the well-formed
-- quasi-smooth tuples from the list of all tuples. The input 'ii' is the index.
-- The result is a list of 7-tuples in the above form. Using the definition of the
-- infinite series, it is easy to write down the corresponding series. Tuples from
-- the tables must also be added, which has not been implemented here.
solve ii = map (filter conds) (makeTuples ii)

lcm3 a b c = lcm a (lcm b c)
gcd3 a b c = gcd a (gcd b c)

-- div   a b gives a/b rounded down
-- divUp a b gives a/b rounded up
divUp a b = -((-a) `div` b)

-- makeTuples - generate all tuples for given index, without checking conditions
makeTuples ii = [makeClass cc ii | cc <- [1..6]]

-- makeClass - generate all tuples for given index and class
makeClass cc ii
  | cc == 1  = concat [makeClassWei a0 (ii - a0)  0 (lcm  a0    (ii - a0)) cc ii |
                       a0 <- [1..(ii `div` 2)]]
  | cc == 2  = concat [makeClassWei a0 a1 (ii - a0) (lcm3 a0 a1 (ii - a0)) cc ii |
                       a0 <- [1..(ii `div` 2)], a1 <- [a0..(ii-a0-1)]]
  | cc == 3  = concat [makeClassWei a0 a1 (ii - a1) (lcm3 a0 a1 (ii - a1)) cc ii |
                       a1 <- [2..(ii `div` 2)], a0 <- [1..(a1-1)]]
  | cc == 4  = concat [makeClassWei (ii-k) (2*k) 0 (lcm (ii-k) k) cc ii |
                       k <- reverse [(max (ii `divUp` 3) 1)..(ii-1)]]
  | cc == 5  = concat [makeClassWei (2*k) (ii-k) 0 (lcm (ii-k) k) cc ii |
                       k <- [1..((ii `divUp` 3)-1)]]
  | cc == 6  = concat [makeClassWei (ii-k) (ii+k) k (lcm3 (ii-k) (ii+k) k) cc ii |
                       k <- reverse [1..(ii-1)]]

-- makeClassWei - create tuples, given smaller weights a0,a1,xx and the number mm
makeClassWei a0 a1 xx mm cc ii
-- for c==1, xx is not used
  | cc == 1  = [Quint a0 a1 b2 b3 (b2+b3) mm cc | b2 <- [a1..(a1+mm-1)],
                b3 <- [b2..(b2+mm-1)]]
-- for c==2, xx = a2
  | cc == 2  = [Quint a0 a1 xx b3 (a1+b3) mm cc | b3 <- [xx..(xx+mm-1)]]
-- for c==3, xx = a2
  | cc == 3  = [Quint a0 a1 xx b3 (a0+b3) mm cc | b3 <- [xx..(xx+mm-1)]]
-- for c==4, xx is not used
  | cc == 4  = [Quint a0 a1 b2 (b2 + (a1 `div` 2)) (2*(b2 + (a1 `div` 2))) mm cc |
                b2 <- [a1..(a1+mm-1)]]
-- for c==5, xx is not used
  | cc == 5  = [Quint a0 a1 b2 (b2 + (a0 `div` 2)) (2*(b2 + (a0 `div` 2))) mm cc |
                b2 <- [a1..(a1+mm-1)]]
-- for c==6, xx = k
  | cc == 6  = [Quint a0 a1 b2 (b2+xx) (a1 + 2*b2) mm cc | b2 <- [a1..(a1+mm-1)]]

-- conds - checks all the well-formedness, quasi-smoothness conditions for a tuple
conds tuple = and[cond j tuple | j <- [1..4]]

-- cond j - checks condition (j) for a given tuple
cond j (Quint a0 a1 a2 a3 dd _ _)
  | j == 1  = and[dd `mod` (gcd ai aj) == 0 | (ai,aj) <- pairs]
  | j == 2  = and[gcd3 ai aj ak == 1 | (ai,aj,ak) <- triples]
  | j == 3  = a3 < dd
  | j == 4  = and[or[(dd - aj) `mod` ai == 0 | aj <- weights] | ai <- weights]
  where
  weights = [a0,a1,a2,a3]
  pairs = [(a0,a1),(a0,a2),(a0,a3),(a1,a2),(a1,a3),(a2,a3)]
  triples = [(a0,a1,a2),(a0,a1,a3),(a0,a2,a3),(a1,a2,a3)]
\end{lstlisting}

\end{document}